\documentclass[]{article}
\begin{document}
\newtheorem{proposition}{Proposition}[section]
\newtheorem{definition}{Definition}[section]
\newtheorem{lemma}{Lemma}[section]

\title{\bf New Characterizations of Algebraic Regularity }
\author{Keqin Liu\\Department of Mathematics\\The University of British Columbia\\Vancouver, BC\\
Canada, V6T 1Z2}
\date{May, 2017}
\maketitle

\begin{abstract} In this paper,  we give new characterizations of algebraic regularity  by using  differential forms and difference quotients.
\end{abstract}

\medskip
\section{Introduction}

In \cite{L}, we introduced the  algebraic regularity on the  quaternions by using a new generalization of Cauchy-Riemann system, and characterized the new generalization of Cauchy-Riemann system by using Fueter-operators. In this paper, we
characterize the  algebraic regularity  on the  quaternions  by using  differential forms and difference quotients.

\medskip
Throughout this paper, we let  $\mathcal{H}=\mathcal{R}e_1\oplus \mathcal{R}e_2\oplus\mathcal{R}e_3\oplus\mathcal{R}e_4$ be the quaternions discovered by W. R. Hamilton in 1843, where $\mathcal{R}$ is the real number field, $e_1$ is the identity of the real division associative algebra $\mathcal{H}$, and the multiplication among the remaining three elements in the $\mathcal{R}$-basis
$\{e_1,\, e_2,\, e_3, \,e_4\}$ is defined by
$$ e_2^2=e_3^2=e_4^2=-e_1, \qquad e_ie_j=-e_je_i=(-1)^{i+j+1}e_{9-i-j},$$
where $2\leq i<j\leq 4$.

\medskip
Recall that the algebraic regularity on the  quaternions  is defined in the following way:

\medskip
\begin{definition}\label{def2.1} Let $U$ be an open subset of $\mathcal{H}$. We say that a quaternion-valued function $f: U\to \mathcal{H}$ is {\bf algebraic regular} at
$c=\displaystyle\sum_{i=1}^4 c_ie_i\in U$ if $f$ has two properties given below.
\begin{description}
\item[(i)] There exist two $C^1$ real-valued functions $f_0: \mathcal{R}^4\to \mathcal{R}$ and $f_1: \mathcal{R}^4\to \mathcal{R}$ such that
\begin{equation}\label{eq1}
f(x)=f_1(x_1,\, x_2,\, x_3,\, x_4)e_1+
\displaystyle\sum_{k=2}^4 x_kf_0(x_1,\, x_2,\, x_3,\, x_4)e_k
\end{equation}
for all $x=\displaystyle\sum_{i=1}^4 x_ie_i\in U$.
\item[(ii)] The following equations hold at $(c_1,\, c_2,\, c_3, \,c_4)\in \mathcal{R}^4$:
\begin{equation}\label{eq2}
\frac{\partial f_1}{\partial x_1}=f_0+x_2\frac{\partial f_0}{\partial x_2}+
x_3\frac{\partial f_0}{\partial x_3}+x_4\frac{\partial f_0}{\partial x_4},
\end{equation}
\begin{equation}\label{eq3}
\frac{\partial f_1}{\partial x_i}=-x_i\frac{\partial f_0}{\partial x_1}, \qquad\qquad
x_i\frac{\partial f_0}{\partial x_j}=x_j\frac{\partial f_0}{\partial x_i},
\end{equation}
where $2\le i,\, j\le 4$ and $i\ne j$.

We say that  $f: U\to \mathcal{H}$ is an {\bf algebraic regular function} on $U$ if $f$ is algebraic regular at every point of $U$.
\end{description}
\end{definition}

\medskip
\section{Characterizing  Algebraic Regularity by Differential Forms}

Let $U$ be an open subset of the quaternion $\mathcal{H}$. A quaternion-valued function $f: U\to \mathcal{H}$ has the form
\begin{equation}\label{eq4}
f(x)=f\left(\displaystyle\sum_{k=1}^4x_ie_k\right)=
\displaystyle\sum_{k=1}^4 f_k(x_1,\, x_2,\, x_3,\, x_4)e_k,
\end{equation}
where $x=\displaystyle\sum_{k=1}^4x_ke_k\in U$, $x_1,\, x_2,\, x_3,\, x_4\in \mathcal{R}$ and $ f_k(x_1,\, x_2,\, x_3,\, x_4)$ is a real-valued function of four real variables $x_1$, $x_2$, $x_3$ and $x_4$. The  quaternion-valued function $f: U\to \mathcal{H}$ given by (\ref{eq4}) is said to be {\bf smooth} ( or $C^1$) if the real-valued function
$ f_k(x_1,\, x_2,\, x_3,\, x_4)$ is smooth (or $C^1$) for $1\le k\le 4$. The alternate notations for the real-value function
$f_k(x_1,\, x_2,\, x_3,\, x_4)$ are given as follows
$$f_k(x_1,\, x_2,\, x_3,\, x_4)=f_k(x)=f_k(x_1e_1+x_ie_i+x_je_j+x_{9-i-j}e_{9-i-j}),$$
where $x=\displaystyle\sum_{k=1}^4x_ke_k$ and $2\le i\ne j\le 4$.

\medskip
For $1\le i\le 4$, we define $dx^i\in Hom_{\mathcal{R}}(\mathcal{H},\,\mathcal{R})$ by
$$(dx^i)(q_1e_1+q_2e_2+q_3e_3+q_4e_4):=q_i,$$
where $q_i\in\mathcal{R}$ for $1\le i\le 4$. Clearly, $\{dx^1,\, dx^2,\, dx^3,\, dx^4\}$ is a basis for the real vector space
$Hom_{\mathcal{R}}(\mathcal{H},\,\mathcal{R})$.

\medskip
A {\bf quaternion-valued $m$-form} on an open subset $U$ of $\mathcal{H}$ is an expression of the form
\begin{equation}\label{eq5}
\alpha=\displaystyle\sum_{1\le i_1, \dots, i_m\le 4}f_{i_1, \dots, i_m}(x)\, dx^{i_1}\wedge  dx^{i_2}\wedge\cdots \wedge
dx^{i_m}\,,
\end{equation}
where $f_{i_1, \dots, i_m}(x): U\to \mathcal{H}$ is a smooth quaternion-valued function on $U$, and
$dx^{i_1}\wedge  dx^{i_2}\wedge\cdots \wedge dx^{i_m}$ is the ordinary  alternating $\mathcal{R}$-multilinear map from
$\underbrace{\mathcal{H}\times\cdots \times\mathcal{H}}_{m}$ to $\mathcal{R}$.   By (\ref{eq5}),  a {\bf quaternion-valued $0$-form} on an open subset $U$ of $\mathcal{H}$ is a  quaternion-valued function defined
on $U$. The smooth quaternion-valued functions $f_{i_1, \dots, i_m}(x)$ are called the
{\bf coefficients} of $\alpha$. We say that $\alpha$ is an {\bf real-valued $m$ form} if all of its coefficients are smooth
real-valued functions.

The {\bf exterior product} of a quaternion-valued $m$-form $\alpha $ given by (\ref{eq5}) and a quaternion-valued $n$-form
$\beta $
\begin{equation}\label{eq7}
\beta=\displaystyle\sum_{1\le j_1, \dots, j_n\le 4}g_{j_1, \dots, j_n}(x)\, dx^{j_1}\wedge  dx^{j_2}\wedge\cdots \wedge
dx^{j_n}
\end{equation}
is defined to be the  quaternion-valued $(m+n)$-form $\alpha\wedge \beta$ which is given by :
$$
\displaystyle\sum_{\begin{array}{c}1\le i_1, \dots, i_m\le 4\\1\le j_1, \dots, j_n\le 4\end{array}}
f_{i_1, \dots, i_m}(x)g_{j_1, \dots, j_n}(x)\, dx^{i_1}\wedge\cdots \wedge
dx^{i_m}\wedge dx^{j_1}\wedge\cdots \wedge  dx^{j_n}.
$$

\bigskip
\begin{proposition}\label{pr3.1}  Let $\alpha$ and $\beta$ be quaternion-valued $m$-form $\alpha $ given by (\ref{eq5}) and a quaternion-valued $n$-form given by (\ref{eq7}). If either  $\alpha$ or $\beta$ is a real-valud form, then
$\alpha\wedge \beta=(-1)^{mn} \beta\wedge \alpha$.
\end{proposition}

\hfill\raisebox{1mm}{\framebox[2mm]{}}

\medskip
If $f: U\to U$ is a quaternion-valued $0$-form given by $f(x)=\displaystyle\sum _{i=1}^4 f_i(x)\,e_i$ with the real-valued functions $f_1(x)$, $f_2(x)$,$f_3(x)$ and $f_4(x)$, we define the {\bf differential } $df$ of $f$ to be the quaternion-valued $1$-form on $U$ given by
\begin{equation}\label{eq11}
df=\frac{\partial f}{\partial x_1}\,dx^1+\frac{\partial f}{\partial x_2}\,dx^2+
\frac{\partial f}{\partial x_3}\,dx^3+\frac{\partial f}{\partial x_4}\,dx^4 ,
\end{equation}
where $\displaystyle\frac{\partial f}{\partial x_i}:=\displaystyle\sum _{j=1}^4
\displaystyle\frac{\partial f_j}{\partial x_i}\,e_j$ is a  quaternion-valued function for $1\le i\le 4$.

\bigskip
\begin{proposition}\label{pr3.2}  If $f$ and $g$ are  quaternion-valued $0$-forms on an open subset $U$ of $\mathcal{H}$, then
\begin{equation}\label{eq12}
d(f+g)=df+dg ,
\end{equation}
\begin{equation}\label{eq13}
d(fg)=f\wedge\, (dg)+(df) \wedge \,g ,
\end{equation}
\begin{equation}\label{eq14}
d(qf)=q\,(df), \qquad d(fq)=(df)\,q \qquad\mbox{for $q\in \mathcal{H}$.}
\end{equation}
\end{proposition}

\hfill\raisebox{1mm}{\framebox[2mm]{}}

\medskip
If $\alpha$ is a quaternion-valued $m$-form given by (\ref{eq5}), we define $d\alpha$ is the  quaternion-valued $(m+1)$-form given by
\begin{equation}\label{eq18}
d\alpha :=\displaystyle\sum_{1\le i_1, \dots, i_m\le 4}d(f_{i_1, \dots, i_m})\wedge
dx^{i_1}\wedge  dx^{i_2}\wedge\cdots \wedge  dx^{i_m}\, .
\end{equation}
The operator $d$ is called the {\bf exterior differentiation}. The next proposition gives the basic properties of the exterior differentiation.

\medskip
\begin{proposition}\label{pr3.3}
\begin{description}
\item[(i)] If $q\in \mathcal{H}$ and  $\alpha$, $\beta$ are  quaternion-valued $m$-forms, then
$d(\alpha +\beta)= d(\alpha)+d(\beta)$, $d(q\alpha)=q\, d(\alpha)$ and $d(\alpha q)=d(\alpha)\, q$.
\item[(ii)] If $\alpha$ is a  quaternion-valued $m$-form and $\beta$ is a  quaternion-valued $n$-form, then
$d(\alpha \wedge \beta)=(d\alpha )\wedge \beta+(-1)^{m}\alpha \wedge (d \beta)$.
\item[(iii)] If $\alpha$ is a  quaternion-valued $m$-form, then $d(d\alpha)=0$.
\end{description}
\end{proposition}

\hfill\raisebox{1mm}{\framebox[2mm]{}}

\bigskip
Following \cite{S}, we use $D\,q$ to denote the following  quaternion-valued $3$-form:
$$
D\,q:=e_1\,  dx^2\wedge  dx^3\wedge  dx^4-e_2\,  dx^1\wedge  dx^3\wedge  dx^4+
e_3\,  dx^1\wedge  dx^2\wedge  dx^4-e_4\,  dx^1\wedge  dx^2\wedge  dx^3
$$
or
\begin{equation}\label{eq44}
D\,q:=e_1\,  dx^2\wedge  dx^3\wedge  dx^4
+\sum_{2\le i<j\le 4}(-1)^{i+j}e_{9-i-j}\,  dx^1\wedge  dx^i\wedge  dx^j .
\end{equation}

Also, recall  from \cite{S} that  left Fueter operator operators $\mathcal{D}_{\ell}:=\displaystyle\sum_{i=1}^4e_i\Big(\frac{\partial  }{\partial x_i}\Big)$ and the right Fueter operator
 $\mathcal{D}_{r}:=\displaystyle\sum_{i=1}^4\Big(\frac{\partial  }{\partial x_i}\Big)e_i$ are defined by
\begin{eqnarray*}
\mathcal{D}_{\ell}(f):&=&e_1\Big(\frac{\partial f }{\partial x_1}\Big)+e_2\Big(\frac{\partial f }{\partial x_2}\Big)
+e_3\big(\frac{\partial f }{\partial x_3}\Big)+e_4\big(\frac{\partial f }{\partial x_4}\Big),\\
\mathcal{D}_r(f):&=&\Big(\frac{\partial f }{\partial x_1}\Big)e_1+\Big(\frac{\partial f }{\partial x_2}\Big)e_2
+\big(\frac{\partial f }{\partial x_3}\Big)e_3+\big(\frac{\partial f }{\partial x_4}\Big)e_4 ,
\end{eqnarray*}
where
$$
\frac{\partial f }{\partial x_i}:=\Big(\frac{\partial f_1 }{\partial x_i}\Big)e_1
+\Big(\frac{\partial f_2 }{\partial x_i}\Big)e_2+
\Big(\frac{\partial f_3 }{\partial x_i}\Big)e_3+\Big(\frac{\partial f_4 }{\partial x_i}\Big)e_4\quad\mbox{for $1\le i\le 4$}.
$$

\medskip
\begin{proposition}\label{pr5.1} If $f: U\to \mathcal{H}$ is a $\mathcal{C}^1$ quaternion-valued function defined on an  open subset of $\mathcal{H}$, then
\begin{equation}\label{eq45}
D\,q\wedge df=-\mathcal{D}_{\ell}(f)\,v\quad\mbox{and}\quad df \wedge D\,q =\mathcal{D}_r(f)\,v ,
\end{equation}
where $v= dx^1\wedge  dx^2\wedge  dx^3\wedge  dx^4$ is the {\bf volume form}.
\end{proposition}

\hfill\raisebox{1mm}{\framebox[2mm]{}}

\bigskip
We now introduce two more real-valued $3$-forms
$D_0\,q$ and  $D_1\,q$  as follows:
\begin{equation}\label{eq47}
D_0\,q:=(x_2+x_3+x_4)\,  dx^2\wedge  dx^3\wedge  dx^4
-\sum_{2\le i<j\le 4}(-1)^{i+j}x_{9-i-j}\,  dx^1\wedge  dx^i\wedge  dx^j ;
\end{equation}
\begin{equation}\label{eq48}
D_1\, q:=  dx^2\wedge  dx^3\wedge  dx^4
+\sum_{2\le i<j\le 4}(-1)^{i+j}\,  dx^1\wedge  dx^i\wedge  dx^j ;
\end{equation}

\medskip
The basic properties of the two  real-valued $3$-forms above   are given in the following
\medskip
\begin{proposition}\label{pr5.2} Let  $U$ be an  open subset of $\mathcal{H}$. If $f$ is a $\mathcal{C}^1$ quaternion-valued function defined on $U$,  then the following equations hold on $U$:
\begin{eqnarray}\label{eq51}
&&D_0\,q\wedge df=- df \wedge D_0\,q \nonumber\\
&=&\Big(-(x_2+x_3+x_4)\frac{\partial f}{\partial x_1}+x_2\frac{\partial f}{\partial x_2}
+x_3\frac{\partial f}{\partial x_3}+x_4\frac{\partial f}{\partial x_4}\Big)\,v,
\end{eqnarray}
\begin{equation}\label{eq52}
D_1\,q\wedge df=- df \wedge D_1\,q =-\Big(\frac{\partial f}{\partial x_1}+
\frac{\partial f}{\partial x_2}+\frac{\partial f}{\partial x_3}+\frac{\partial f}{\partial x_4}\Big)\,v,
\end{equation}
where $v= dx^1\wedge  dx^2\wedge  dx^3\wedge  dx^4$ is the {\bf  volume form}.
\end{proposition}

\hfill\raisebox{1mm}{\framebox[2mm]{}}

\bigskip
Using Proposition \ref{pr5.2} , we have

\medskip
\begin{proposition}\label{pr5.3} Let  $U$ be an open subset of $\mathcal{H}$. If $f: U\to \mathcal{H}$ is a function given by
$$f(x)=f_1(x_1, x_2, x_3, x_4)e_1+\displaystyle\sum_{k=2}^4 x_kf_0(x_1, x_2, x_3, x_4)e_k ,$$
where $x=\displaystyle\sum_{i=1}^4 x_ie_i\in U$ with $x_1$, $x_2$, $x_3$, $x_4 \in \mathcal{R}$, $f_1$ and $f_0$ are
$\bf{C}^1 $ functions, then the following are equivalent:
\begin{description}
\item[(i)] $f$ is algebraic regular on $U$;
\item[(ii)] Both the equation
\begin{equation}\label{eq76}
D\,q\wedge df+2\, D_0\,q\wedge df_0+2\, D_1\,q\wedge df_1=0
\end{equation}
and the the equation
\begin{equation}\label{eq77}
df\wedge D\,q+2\, df_0\wedge D_0\,q+2\, df_1\wedge D_1\,q=0
\end{equation}
hold on $U$.
\end{description}
\end{proposition}

\hfill\raisebox{1mm}{\framebox[2mm]{}}

\medskip
\section{Characterizing  Algebraic Regularity by  Difference Quotients}

Let $f(x)=f_1(x_1, x_2, x_3, x_4)e_1+\displaystyle\sum_{k=2}^4 x_kf_0(x_1, x_2, x_3, x_4)e_k $ be a quaternion-valued function defined on an open subset $U$ of $\mathcal{H}$, where $x=\displaystyle\sum_{i=1}^4 x_ie_i$ with $x_1$, $x_2$, $x_3$, $x_4\in \mathcal{R}$ for
$1\le i\le 4$. For each $c=\displaystyle\sum_{i=1}^4 c_ie_i$ with $c_1$, $c_2$, $c_3$, $c_4\in \mathcal{R}$, we define six pure
quaternion-valued functions on $U$ as follows:
\begin{eqnarray}\label{eq113}
\stackrel{\leftarrow}{f}^c_{9-i-j}(x):& =&\big[c_i\,f_0(c_1e_1+x_ie_i+c_je_j+c_{9-i-j}e_{9-i-j})+\nonumber\\
&&\quad +c_j\,f_0(c_1e_1+c_ie_i+x_je_j+c_{9-i-j}e_{9-i-j})\big]e_{9-i-j}+\nonumber\\
&&  -\big[(-1)^{i+j}c_if_0(x_1, c_2, c_3, c_4)+\nonumber\\
&&\quad +c_{9-i-j}\,f_0(c_1e_1+c_ie_i+x_je_j+c_{9-i-j}e_{9-i-j})\big]e_j+\nonumber\\
&&  +\big[(-1)^{i+j}c_jf_0(x_1, c_2, c_3, c_4)+\nonumber\\
&&\quad  -c_{9-i-j}\,f_0(c_1e_1+x_ie_i+c_je_j+c_{9-i-j}e_{9-i-j})\big]e_i ,
\end{eqnarray}
\begin{eqnarray}\label{eq114}
\stackrel{\rightarrow}{f}^c_{9-i-j}(x): & =&\big[c_i\,f_0(c_1e_1+x_ie_i+c_je_j+c_{9-i-j}e_{9-i-j})+\nonumber\\
&&\quad +c_j\,f_0(c_1e_1+c_ie_i+x_je_j+c_{9-i-j}e_{9-i-j})\big]e_{9-i-j}+\nonumber\\
&&  +\big[(-1)^{i+j}c_if_0(x_1, c_2, c_3, c_4)+\nonumber\\
&&\quad -c_{9-i-j}\,f_0(c_1e_1+c_ie_i+x_je_j+c_{9-i-j}e_{9-i-j})\big]e_j+\nonumber\\
&&  -\big[(-1)^{i+j}c_jf_0(x_1, c_2, c_3, c_4)+\nonumber\\
&&\quad  +c_{9-i-j}\,f_0(c_1e_1+x_ie_i+c_je_j+c_{9-i-j}e_{9-i-j})\big]e_i,
\end{eqnarray}
where $2\le i<j\le 4$.

\medskip
In the proposition below, we characterize  the algebraic regularity by the limits of a new kind of difference quotients which use the six pure quaternion-valued functions $\stackrel{\leftarrow}{f}^c_i(x)$ and $\stackrel{\rightarrow}{f}^c_i(x)$ with
$i=2$, $3$ and $4$.

\medskip
\begin{proposition}\label{pr6.1} Let  $f: U\to \mathcal{H}$ be a quaternion-valued function defined by
$$f(x)=f_1(x_1, x_2, x_3, x_4)e_1+\displaystyle\sum_{i=2}^4 x_if_0(x_1, x_2, x_3, x_4)e_i ,$$
where $U$ is an open subset of $\mathcal{H}$, $x=\displaystyle\sum_{k=1}^4x_ke_k\in U$ with  $x_1$, $x_2$, $x_3$, $x_4\in \mathcal{R}$. Let
$c=\displaystyle\sum_{k=1}^4c_ke_k\in U$ and $\Delta q=\displaystyle\sum_{k=1}^4(\Delta q)_ke_k$, where $c_k$,
$(\Delta q)_k\in \mathcal{R}$ for $1\le k\le 4$. If  $f_1(x_1, x_2, x_3, x_4)$ and $f_0(x_1, x_2, x_3, x_4)$ are
$\mathcal{C}^1 $ functions, then the following are equivalent:
\begin{description}
\item[(i)] $f$ is algebraic regular at $x=c$;
\item[(ii)] $\displaystyle\lim_{\Delta q\to 0}(\Delta q)^{-1}\Big\{f(c+\Delta q)-f(c)+
\displaystyle\sum_{i=2}^4\Big[\stackrel{\leftarrow}{f}^c_i\Big(c+(\Delta q)_i\displaystyle\sum_{k=1}^4e_k\Big)-
\stackrel{\leftarrow}{f}^c_i(c)\Big]\Big\}$ exists;
\item[(iii)] $\displaystyle\lim_{\Delta q\to 0}\Big\{f(c+\Delta q)-f(c)+
\displaystyle\sum_{i=2}^4\Big[\stackrel{\rightarrow}{f}^c_i\Big(c+(\Delta q)_i\displaystyle\sum_{k=1}^4e_k\Big)-
\stackrel{\rightarrow}{f}^c_i(c)\Big]\Big\}(\Delta q)^{-1}$ exists.
\end{description}
Moreover, if one of the three coditions above holds, then both the limit in (ii) and the limit in (iii) equal to
$\displaystyle\frac{\partial f}{\partial x_1}(c)$.
\end{proposition}

\hfill\raisebox{1mm}{\framebox[2mm]{}}

\bigskip
Based on Proposition \ref{pr6.1}, we call $\displaystyle\frac{\partial f}{\partial x_1}$ the {\bf quaternion derivative} of an algebraic regular function $f$ on an open subset $U$ of $\mathcal{H}$. It is easy to check that if  $f$ is an algebraic regular function on an open subset $U$  of $\mathcal{H}$, then its quaternion derivative
$\displaystyle\frac{\partial f}{\partial x_1}$  is also an algebraic regular function on the open subset $U$.


\bigskip
\begin{thebibliography}{9}
\bibitem{F} R. Fueter, \textsl{Die Funktionentheorie der Differentialeichungen $\Delta u=0$ and
$\Delta\Delta u=0$ mit vier reellen Variablen}, Comment. Math. Helv. 7 (1935), 307-330
\bibitem{L} K. Liu, \textsl{Algebraic Regularity over Quaternions and Regular Four-Manifolds}, arXiv:1511.08532
\bibitem{S} A. Sudbery, \textsl{Quaternionic analysis}, Math. Proc. Camb. Phil. Soc. 85 (1979), 199-225
\end{thebibliography}
\end{document}